\documentstyle{amsppt}
\magnification1200
\pagewidth{6.5 true in}
\pageheight{9.25 true in}
\NoBlackBoxes


\topmatter
\title Degree $1$ elements of the Selberg class
\endtitle
\author K. Soundararajan
\endauthor
\thanks{The author is partially supported by the American Institute of
Mathematics and by the
National Science Foundation.}\endthanks
\address
Department of Mathematics,  University of Michigan,
Ann Arbor, MI 48109, USA
\endaddress
\email
ksound\@umich.edu
\endemail
\endtopmatter

\def\phi{\varphi}
\def\lam{\lambda}

\document

\noindent   In [5] A. Selberg axiomatized properties expected of
$L$-functions and introduced the ``Selberg class.'' 
We recall that an element $F$ of the Selberg class ${\Cal S}$ satisfies the
following axioms.

\noindent {\bf Axiom 1}. In the half-plane $\sigma>1$ the
function $F(s)$ is given by an absolutely convergent Dirichlet
series $\sum_{n=1}^{\infty} a(n) n^{-s}$ with $a(1)=1$ and
$a(n)\ll n^{\epsilon}$ for every $\epsilon >0$.

\noindent {\bf Axiom 2}. There is a natural number $m$
such that $(s-1)^{m} F(s)$ extends to an analytic function in the
entire complex plane.

\noindent {\bf Axiom 3}. There is a function
$\Phi(s) = Q^{s} G(s) F(s)$ where $Q >0$ and
$$
G(s) = \prod_{j=1}^{r} \Gamma(\lam_j s+\mu_j)
\qquad \text{with} \qquad \lam_j>0 \text{  and  } \text{Re }\mu_j \ge 0
$$
such that
$$
\Phi(s) = \omega \overline{\Phi}(1-s),
$$
where $|\omega| =1$ and for any function $f$ we denote $\overline{f}(s)
= \overline{f(\overline{s})}$. We let $d:= 2\sum_{j=1}^{r} \lam_j$ denote
the ``degree'' of $F$.

\noindent {\bf Axiom 4}. We may express $\log F(s)$ by a Dirichlet series
$$
\log F(s) = \sum_{n=2}^{\infty} \frac{b(n)}{n^s} \frac{\Lambda(n)}{\log n}
$$
where $b(n) \ll  n^{\vartheta}$ for some $\vartheta <\frac 12$.  Set 
$b(n)=0$ if $n$ is not a prime power. 


A fundamental conjecture asserts that the degree of an element in
the Selberg class is an integer.  From the work of H.E. Richert [4] it
follows that there are no elements in the Selberg class with
degree $0< d<1$.  This was rediscovered by J.B. Conrey and A. Ghosh [1]
who also proved that the Selberg class of degree $0$ contains only the
constant function $1$.  Recently J. Kaczorowski and A. Perelli [2]
determined the structure of the Selberg class for degree $1$ and
showed that this consists of the Riemann zeta function and shifts of
Dirichlet $L$-functions.  Subsequently in [3] they showed that there
are no elements of the Selberg class with degree $1< d <5/3$.  In this
note we shall give a short and simple proof of Kaczorowski and
Perelli's beautiful result on the Selberg class for degree $1$.

\proclaim{Theorem}  Suppose $F$ satisfies Axioms 1 to 3 and that the
degree of $F$ is $1$.  Then there exists a positive integer $q$ and
a real number $A$ such that
$a(n) n^{-iA}$ is periodic $\pmod q$.  If in addition $F$ satisfies
Axiom 4 then there is a primitive Dirichlet character $\chi'
\pmod {q'}$ such that $F(s) = L(s+iA,\chi')$.
\endproclaim

We remark that Kaczorowski and Perelli obtain their results 
without assuming the hypothesis $a(n) \ll n^{\epsilon}$.  
We could restructure our proof to avoid this assumption, but have 
preferred not to do so in the interest of keeping the exposition 
transparent.  Our method may also be modified and combined with 
the ideas in Kaczorowski and Perelli [3] to give a simplification 
of their result for $1<d<5/3$.


Suppose $F$ satisfies Axioms 1 to 3 and has degree $1$.  By
Stirling's formula we see that for $t \ge 1$
$$
\frac{\overline{G}(1/2-it)}{G(1/2+it)} = e^{-it \log \frac{t}{2e}
+i \frac{\pi}{4} +iB} t^{iA} C^{-it} \Big( 1+ O\Big(\frac 1t\Big)\Big),
\tag{1}
$$
for some real numbers $A$, $B$ and $C>0$.  Let $\alpha$ be positive and
$T\ge 1$.  Define
$$
{\Cal F}(\alpha,T) = \frac{1}{\sqrt{\alpha}} \int_{\alpha
T}^{2\alpha T} F(1/2+it) e^{it \log \frac{t}{2\pi e\alpha}
-i\frac{\pi}{4} } dt, \tag{2a}
$$
and set (it will follow from our proof that the limit below is
well defined)
$$
{\Cal F}(\alpha) = \lim_{T \to \infty} \frac{1}{T^{1+iA}} {\Cal
F}(\alpha, T). \tag{2b}
$$

\proclaim{Lemma}  For any real number $t$ and all $X\ge 1$ we have
$$
F(\tfrac 12+it) = \sum_{n=1}^{\infty} \frac{a(n)}{n^{\frac 12+it}}
e^{-n/X} + O((1+|t|)^{1+\epsilon} X^{-1+\epsilon} + X^{\frac
12+\epsilon}e^{-|t|}).
$$
\endproclaim
\demo{Proof}  Consider for $c>\tfrac 12$
$$
\frac{1}{2\pi i} \int_{c-i\infty}^{c+i\infty} F(\tfrac 12+it
+w)X^w \Gamma(w) dw.
$$
Expanding $F(\tfrac 12+it +w)$ into its Dirichlet series and integrating
term by term we see
that this equals $\sum_{n} a(n) n^{-\frac 12-it} e^{-n/X}$. Next
we move the line of integration to Re$(w)= -1+\epsilon$.  The pole at
$w=0$ leaves the residue $F(\frac 12+it)$.  The possible pole at
$w=\frac 12 -it$ leaves a residue $\ll X^{\frac 12+ \epsilon}
(1+|t|)^{\epsilon} |\Gamma(\frac 12+it)| \ll X^{\frac 12+\epsilon}
e^{-|t|}$ due to the rapid decay of $\Gamma (\frac 12+it)$. Note
that by the functional equation and Stirling's formula $|F(\frac
12 +it+w)| = Q^{1-2\text{Re }w} \frac{|\overline{G}(\tfrac 12
-it-w)|}{|G(\tfrac 12+it+w)|} |F(\tfrac 12-it -w)| \ll
(1+|t|+|w|)^{1+\epsilon}$ for any $w$ on the line Re
$w=-1+\epsilon$.  Hence the integral on the line Re$(w)=-1+\epsilon$
is $ \ll X^{-1+\epsilon} (1+|t|)^{1+\epsilon}$.
\enddemo

Using the functional equation in (2a) we see that
$$
{\Cal F}(\alpha, T) = \frac{\omega}{\sqrt{\alpha} }
\int_{\alpha T} ^{2\alpha T} {\overline F}(1/2-it) Q^{-2it}
\frac{\overline{G}(1/2-it)}{G(1/2+it)} e^{it\log \frac{t}{2\pi
e\alpha}-i\frac{\pi}{4}} dt
$$
and using (1) this is
$$
 \frac{\omega e^{iB}}{\sqrt{\alpha}}
\int_{\alpha T}^{2\alpha T} {\overline F}(1/2-it) (\pi C Q^2
\alpha)^{-it} t^{iA}(1+O(1/T)) dt.
$$
We now input our Lemma above with $X=T^{\frac 43}$ to deduce that
$$
\align
 {\Cal F}(\alpha, T) &= \frac{\omega e^{iB}}{\sqrt{\alpha}}
\int_{\alpha T}^{2\alpha T} \sum_{m}
\frac{\overline{a(m)}}{\sqrt{m}} e^{-m/X} \Big(\frac{m}{\pi C
Q^2\alpha}\Big)^{it} t^{iA}
\Big( 1+O\Big(\frac 1T\Big)\Big) dt +O(T^{\frac 23 +\epsilon})\\
&=\frac{\omega e^{iB}}{\sqrt{\alpha}} \sum_{m}
\frac{\overline{a(m)}}{\sqrt{m}} e^{-m/X} \int_{\alpha T}^{2\alpha
T} \Big( \frac{m}{\pi C Q^2\alpha}\Big)^{it} t^{iA} dt + O(T^{\frac
23+\epsilon} ).
\tag{3}
\\
\endalign
$$

If $x\neq 1$ then integration by parts gives that
$$
\int_{\alpha T}^{2\alpha T} x^{it} t^{iA} dt
= \frac{(2\alpha T)^{iA} x^{2iT} - (\alpha T)^{iA}
x^{iT}}{i\log x} -\int_{\alpha T}^{2\alpha T}
\frac{x^{it}}{i\log x} iA t^{iA -1}
dt \ll \frac{1}{|\log x|}, \tag{4a}
$$
while if $x=1$ we have that
$$
\int_{\alpha T}^{2\alpha T} t^{iA} dt = \frac{(2\alpha T)^{1+iA}-
(\alpha T)^{1+iA}}{1+iA}.
\tag{4b}
$$
Using (4a,b) in (3) we obtain that
$$
\align
{\Cal F}(\alpha)
&= \lim_{T\to \infty} T^{-1-iA}{\Cal F}(\alpha, T) \\
&= \omega
e^{iB} \delta(\pi C Q^2 \alpha \in {\Bbb N})
\frac{\overline{a(\pi C Q^2\alpha)}\alpha^{iA}}{\sqrt{\pi C} Q}
\frac{2^{1+iA}
-1}{1+iA}+O\Big( \lim_{T\to \infty} T^{-1+\epsilon} \sum_{m}
\frac{|a(m)|}{\sqrt{m}} e^{-m/X}\Big) \\
&= \omega
e^{iB} \delta(\pi CQ^2 \alpha \in {\Bbb N})
\frac{\overline{a(\pi CQ^2\alpha)}\alpha^{iA}}{\sqrt{\pi C} Q} \frac{2^{1+iA}
-1}{1+iA},
\tag{5}\\
\endalign
$$
where $\delta(\pi C Q^2\alpha \in {\Bbb N}) =1$ if $\pi C Q^2
\alpha \in {\Bbb N}$ and is $0$ otherwise.

We now present a different way of evaluating ${\Cal F}(\alpha, T)$
which will show that ${\Cal F}(\alpha)$ is periodic in $\alpha$ with
period $1$.  Using our Lemma with $X= T^{\frac 43}$ again, we see that
$$
{\Cal F}(\alpha,T)
= \frac{1}
{\sqrt{\alpha}} \sum_{n} \frac{a(n)}{\sqrt{n}}
e^{-n/X} \int_{\alpha T}^{2\alpha T} e^{it\log \frac{t}{2\pi e n\alpha}
-i\frac {\pi}{4} }dt
+ O(T^{\frac 23+\epsilon}). \tag{6}
$$
The oscillatory integral in (6) above is estimated by familiar
techniques, see Lemmas 4.2 and 4.6 of E.C. Titchmarsh [6]. For
$2\pi n > 3T$ we use Lemma 4.2 of Titchmarsh which shows
that the integral is $\ll 1$.  Thus the contribution
of such $n$ to (6) is $\ll T^{\frac{2}{3} +\epsilon}$.
For smaller $n$ we use Lemma 4.6
of Titchmarsh.  In the range $T\le 2\pi n \le 2T$ we obtain that the
integral is
$$
2\pi \sqrt{n\alpha} e(-n\alpha) + O\Big( T^{\frac 25} + \min
\Big(\sqrt{T},\frac{1}{|\log (T/2\pi n)|} \Big) + \min
\Big(\sqrt{T}, \frac{1}{|\log (T/\pi n)|}\Big) \Big).
$$
For $2\pi n$ below $T$ or between $2T$ and $3T$ the integral
is bounded by the error terms above.  Piecing this together we
conclude that
$$
{\Cal F}(\alpha, T) = 2\pi \sum_{T\le 2\pi n \le 2T} a(n)
e(-n\alpha) + O(T^{\frac 9{10}+\epsilon}).
$$
Hence ${\Cal F}(\alpha)= {\Cal F}(\alpha+1)$ which leads by (5)
to
$$
\delta(\pi C Q^2\alpha \in {\Bbb N}) \overline{a(\pi C Q^2\alpha)}
\alpha^{iA} = \delta(\pi C Q^2(\alpha+1) \in {\Bbb N})
\overline{a(\pi C Q^2(\alpha+1))}  (\alpha +1)^{iA}. \tag{7}
$$

From (7) we deduce immediately that $\pi C Q^2 =q$ must be a positive
integer and further that $\overline{a(n)} n^{iA}$ must be
periodic $\pmod q$.  This proves the first part of our Theorem.

Suppose now that $F$ also satisfies Axiom 4 so that the
coefficients $a(n)$ are multiplicative.  Periodicity and
multiplicativity together imply that for all $n$ coprime to $q$,
$a(n)n^{-iA}$ must equal $\chi(n)$ for a Dirichlet character $\chi
\pmod q$.  Let $\chi^{\prime} \pmod {q'}$ denote the primitive
character inducing $\chi$. Then the ratio
$F(s)/L(s+iA,\chi^{\prime})$ is an Euler product over the finitely
many primes dividing $q$, and by Axiom (4) the logarithm of this
Euler product converges absolutely in the half plane $\sigma
> {\vartheta}$ (recall that $\vartheta <\frac 12$).  Thus, with ${\frak a} 
= (1-\chi^{\prime}(-1))/2$,  
$$
H(s):= \frac{Q^s G(s) F(s)}{(q'/\pi)^{\frac s2} 
\Gamma(\frac{s+iA+{\frak a}}{2}) 
L(s+iA,\chi')}
$$
is an entire function in Re$(s) > \vartheta$ and $\overline{H}(s)$ is also 
entire in this region.  Further both functions are free of zeros in this 
region.  The functional equations for $F$ and $L$ now show that 
$H(s)$ and $\overline{H}(s)$ are entire and free of zeros in the region 
Re$(s) < 1-\vartheta$.  Since $\vartheta <\frac 12$ we deduce that 
$H(s)$ and ${\overline H}(s)$ are entire functions in all of ${\Bbb C}$ and 
that they never vanish.  Since $H$ is the ratio of two entire functions
\footnote{If $\chi'$ is the trivial character then multiply 
the numerator and denominator in the definition of $H$ by $s(s-1)$ to 
make them regular at $1$.} 
of order $1$ it follows by Hadamard's theorem that $H(s) = 
ae^{bs}$ for some constants $a$ and $b$.  The 
functional equation connecting $H(s)$ and ${\overline{H}}(1-s)$ now mandates 
that $b=0$ and so $H$ is a constant.   Examining the behaviour of 
$H(1/2+it)$ for large $t$ it follows easily that 
$F(s)=L(s+iA,\chi^{\prime})$, proving our Theorem.

\enddocument